\documentclass{amsart}
\usepackage{amsmath,amssymb,latexsym}
\usepackage[latin1]{inputenc}
\usepackage[all]{xy}
\usepackage[dvips]{graphics}

\newtheorem{theorem}{Theorem}[section]
\newtheorem{lemma}[theorem]{Lemma}
\newtheorem{proposition}[theorem]{Proposition}
\newtheorem{corollary}[theorem]{Corollary}

\theoremstyle{definition}
\newtheorem{definition}[theorem]{Definition}
\newtheorem{example}[theorem]{Example}

\theoremstyle{remark}
\newtheorem{remark}[theorem]{Remark}
\newtheorem{conjecture}[theorem]{Conjecture}

\newcommand{\N}{\mathbb N}
\newcommand{\Z}{\mathbb Z}
\newcommand{\Q}{\mathbb Q}
\newcommand{\C}{\mathbb C}
\newcommand{\fq}{{\mathbb F}_q}

\newcommand{\dimv}{\underline{\dim}}

\newcommand{\Gr}{\operatorname{Gr}}
\newcommand{\id}{\operatorname{Id}}

\newcommand{\qc}{\mathcal Q}
\newcommand{\F}{\mathcal F}
\newcommand{\A}{\mathcal A}
\newcommand{\ab}{\A(B)}
\newcommand{\ad}{\A(\Delta)}

\newcommand{\ext}{\operatorname{Ext}^1}
\renewcommand{\hom}{\operatorname{Hom}}

\newcommand{\Rad}{\operatorname{Rad}}
\newcommand{\Ind}{\operatorname{Ind}}
\newcommand{\Mod}{\operatorname{mod}}
\newcommand{\mkq}{\Mod_k(Q)}
\newcommand{\imkq}{\overline{\Mod}_k(Q)}
\newcommand{\indkq}{\Ind_k(Q)}

\newcommand{\ud}{\underline d}
\newcommand{\ue}{\underline e}
\newcommand{\uf}{\underline f}
\newcommand{\ug}{\underline g}
\newcommand{\ux}{\underline x}
\newcommand{\uu}{\underline u}
\newcommand{\um}{\underline m}
\newcommand{\un}{\underline n}

\newcommand{\phap}{\Phi_{\geq -1}}
\newcommand{\edge}{\mathsf{E}}
\newcommand{\source}{\mathsf{S}}
\newcommand{\cld}{\operatorname{Cl}_{\Delta}}

\newcommand{\gremk}{\Gr_{\ue}(M)_k}
\newcommand{\grem}{\Gr_{\ue}(M)}
\newcommand{\gq}{\Gamma_Q}
\newcommand{\gcc}{\Gamma_{\cc}}
\newcommand{\cc}{{\mathcal C}}
\newcommand{\ccd}{{\mathcal C}_\Delta}
\newcommand{\ccq}{{\mathcal C}_Q}

\begin{document}
\title[Cluster algebras as Hall algebras]{Cluster algebras as Hall algebras of quiver representations}
\author{Philippe Caldero}
\address{Institut Camille Jordan, Université Claude Bernard Lyon I,
69622 Villeurbanne Cedex, France}
\email{caldero@igd.univ-lyon1.fr}
\author{Frédéric Chapoton}
\address{Institut Camille Jordan, Université Claude Bernard Lyon I,
69622 Villeurbanne Cedex, France}
\email{chapoton@igd.univ-lyon1.fr}

\begin{abstract}
  Recent articles have shown the connection between representation
  theory of quivers and the theory of cluster algebras.  In this
  article, we prove that some cluster algebras of type ADE can be
  recovered from the data of the corresponding quiver representation
  category. This also provides some explicit formulas for cluster
  variables.
\end{abstract}

\date{\today}

\maketitle
\section{Introduction}
Cluster algebras were introduced in \cite{FZ1} by S. Fomin and A.
Zelevinsky in connection with the theory of dual canonical bases and
total positivity. Coordinate rings of many varieties from Lie group
theory -- semisimple Lie groups, homogeneous spaces, generalized
Grassmannian, double Bruhat cells, Schubert varieties -- have a
structure of cluster algebra, at least conjecturally, see
\cite{BFZ,scott}. One of the goals of the theory is to provide a
general framework for the study of canonical bases of these coordinate
rings and their $q$-deformations.\par
A (coefficient-free) cluster algebra $\A$ of rank $n$ is a subalgebra
of the field $\Q(u_1,\ldots,u_n)$. It is defined from a distinguished
set of generators, called cluster variables, constructed by an
induction process from a antisymmetrizable matrix $B$, see Section
\ref{recoll}.  The Laurent phenomenon asserts that $\A$ is a
subalgebra of $\Q[u_1^{\pm 1},\ldots,u_n^{\pm 1}]$. There exists a
notion of compatibility between two cluster variables; maximal subsets
of pairwise compatible cluster variables are called clusters. All
clusters have the same cardinality, which is the rank of the cluster
algebra.\par
A cluster algebra is of {\it finite type} if the number of cluster
variables is finite. The classification of cluster algebras of finite
type \cite{FZ2} is a fundamental step in the theory. The main result
is that these cluster algebras come from an antisymmetrized Cartan
matrix of finite type, see Section \ref{finite}. Moreover, in this
case the cluster variables are in correspondence with the set of
almost positive roots $\phap$, \textit{i.e.}  positive roots or
opposed simple roots, of the root system.\par
The Gabriel theorem asserts that the set of indecomposable
representations of a quiver $Q$ of Dynkin type is in bijection with
the set $\Phi_+$ of positive roots. The cluster category $\cc$ was
constructed in \cite{BMRRT,CCS} as an extension of the category $\mkq$
of finite dimensional $kQ$-modules, such that the set of
indecomposable objects of $\cc$ is in bijection with $\phap$. The
category $\cc$ is not abelian in general, but it is a triangulated
category, \cite{keller}. In \cite{BMRRT}, this category is studied in
depth. The authors give a correspondence between cluster variables and
indecomposable objects of $\cc$. They prove that the compatibility of
two cluster variables correspond to the vanishing of the Ext groups;
hence, clusters correspond to so-called ext-configurations. They prove
that there exist many analogies between finite cluster algebras and
cluster categories, but the properties of the correspondence are
mostly conjectural, see \cite[Conjecture 9.3]{BMRRT}.\par
In \cite{CCS}, the authors prove that the denominators of cluster
variables can be calculated from $\cc$ in type A. They give a
combinatorial/geometric approach of $\cc$ in the spirit of Teichmüller
spaces, \cite{fock}.\par
The implicit question behind all articles \cite{MRZ,BMRRT,CCS} dealing
with cluster algebras and quiver theory is: can one realize the
cluster algebra as a ``Hall algebra'' of the category $\cc$ in some
sense ?\par
Recall that $\mkq$ is a (non full) subcategory of $\cc$. In this
article, the cluster variable associated to an indecomposable
$kQ$-module, in fact to any module, is explicitly given, see Theorem
\ref{main}. This result is interesting from different
angles.\par\noindent 1. It strengthens the relations between the
category $\cc$ and the algebra $\A$.\par\noindent 2. We obtain here
explicit expressions for cluster variables, instead of inductive ones.
These expressions are in terms of Euler-Poincar\'e characteristic of
Grassmannians of submodules. Note that these characteristics can be
easily calculated in the A$_n$ case, see Example \ref{An}. They can
also be calculated in a combinatorial way in the D$_n$
case.\par\noindent 3. One important open question in cluster theory is
the positivity conjecture \cite[\S 3]{FZ1}, which says that cluster
variables should be Laurent polynomials with positive coefficients in
the variables of any fixed cluster. Our explicit expressions will be
used in another article \cite{caldkell} to show that cluster variables
indeed have a positive Laurent expansion in any cluster associated to
a Dynkin type quiver. Recall that the positivity conjecture is known
to hold only for a distinguished cluster so far \cite[Th. 1.10]{FZ2}.
\par\noindent 
4. The expression gives the possibility to quantize
cluster algebras in the Ringel-Hall algebras spirit: the
Euler-Poincar\'e characteristic should be replaced by a polynomial
which counts $\fq$-rational points on the variety.\par
These points as well as other ones (toric degenerations, denominator theorems) will be developped in a forthcoming article, \cite{caldkell}, which mainly relies on the cluster variable formula, Theorem
\ref{main}.\par 
In the sequel, we give a conjectural expression for cluster variables
associated to a multiplicity-free indecomposable module over any
quiver of simply-laced finite cluster type. As a special case, this
conjecture enables, in A$_n$ type, to calculate in a combinatorical
way the cluster variables in terms of any cluster.\par
To conclude, we give a connection between our theorem, the geometric
realization of \cite{CCS}, and the Coxeter-Conway friezes
\cite{CoxCon}.

\smallskip

We would like to note that one of the starting points for the
experimental work leading to this article was the combinatorial
expressions for some $Y$-system Laurent polynomials given in
\cite{Ysystems} for multiplicity-free roots. Although these are
definitely not the same as cluster Laurent polynomials, the
combinatorics is quite similar.

\smallskip
    
{\bf Acknowledgments}: The first author would like to thank Markus
  Reineke for conversations on Euler-Poincaré characteristic and
  Grassmannians of submodules. He is also grateful to Bernhard Keller
  for a simpler argument in the proof of Lemma \ref{merciK}. 

\section{Recollection from cluster algebras and cluster categories}
\subsection{ }\label{recoll}
In this section, we give basic definitions and Theorems concerning
cluster algebras, see \cite{FZ1,FZ2,BFZ}. The cluster algebras in this
article are defined on a trivial semigroup of coefficients, and will
be called {\it reduced} cluster algebras. The recollection below is
transposed into the framework of reduced cluster algebras.\par
Let $n$ be a positive integer and let $B=(b_{ij})$ be a square matrix
in $M_n(\Z)$. We say that $B$ is antisymmetrizable if there exists a
diagonal matrix in $M_n(\N)$ such that $DB$ is antisymmetric. We
introduce the field $\F:=\Q(u_1,\ldots,u_n)$, with algebraically
independent generating set $\uu:=(u_1,\ldots,u_n)$. A pair $(\ux, B)$,
where $\ux=(x_1,\ldots,x_n)$ is an algebraically independent
generating set of $\F$ and where $B$ is an antisymmetrizable matrix,
will be called a seed. In the sequel, we will identify the rows and
the columns of the matrix $B$ with the elements of $\ux$.\par
Fix a seed $(\ux,B)$, $B=(b_{yz})$, and $w$ in the base $\ux$. Let
$w'$ in $\F$ be such that
\begin{equation}\label{exchange}
ww'=\prod_{b_{yw}>0}y^{b_{yw}}+\prod_{b_{yw}<0}y^{-b_{yw}}.
\end{equation}
This is the so-called exchange relation. Now, set
$\ux':=\ux-\{w\}\cup\{w'\}$ and $B'=(b_{yz}')$ such that
\begin{equation}
  b_{yz}'=\left\{
    \begin{array}{ll}
      -b_{yz} &\textup{if } y=w \textup{ or } z=w,\\
      b_{yz}+1/2(\mid b_{yw}\mid b_{wz}+b_{yw}\mid b_{wz}\mid) &
      \textup{otherwise.}
    \end{array}\right.
\end{equation}
Then, it is known that $(\ux', B')$ is also a seed. We say that this
seed is the mutation of the seed $(\ux,B)$ in the direction $w$. We
also say that $w$ and $w'$ form an exchange pair. It is easily seen
that the mutation of the seed $(\ux',B')$ in the direction $w'$ is
$(\ux,B)$. We can define the equivalence relation generated by
$(\ux,B)\sim(\ux',B')$ if $(\ux',B')$ is a mutation of $(\ux,B)$.\par
We assign to a antisymmetrizable matrix $B$ a $\Q$-algebra in the
following way.
\begin{definition}
  The reduced cluster algebra $\ab$ associated to the
  antisymmetrizable matrix $B$ is the subalgebra of $\F$ generated by
  all $\ux$ such that $(\uu,B)\sim(\ux,B')$. Such $\ux$ are called
  clusters and the elements of $\ux$ are called cluster variables.
\end{definition}
\begin{remark}
  More generally, see \cite{BFZ}, cluster algebras are associated to
  rectangular matrices in $M_{n,m}(\Z)$. We will not be concerned with
  such algebras in this article.
\end{remark}
Note the so-called {\it Laurent phenomenon}, see \cite{FZ1}:
\begin{theorem}
  Let $B$ be a antisymmetrizable matrix in $M_n(\Z)$, then $\ab$ is a
  subalgebra of $\Q[u_i^{\pm 1}, 1\leq i\leq n]$.
\end{theorem}

\subsection{ }\label{finite}
This section is concerned with {\it finite} reduced cluster algebras,
\textit{i.e.} cluster algebras with a finite number of cluster
variables.\par
Let $\Delta$ be a Dynkin diagram of rank $n$ and let A$_\Delta$ be its
Cartan matrix. We denote by $\Phi$, resp. $\Phi_+$, the root system,
resp. the set of positive roots, associated to $\Delta$. Let
$\alpha_i$, $1\leq i\leq n$, be the simple roots and let $\qc$ be the
$\Z$-lattice generated by them. We also denote by $\phap$ the set of
\textit{almost positive} roots
$\Phi_+\cup\{-\alpha_1,\ldots,-\alpha_n\}$.\par
We have the following fundamental Theorem, see \cite{FZ2}:
\begin{theorem}
  A reduced cluster algebra $\A$ is finite if and only if there exists
  a seed $(\ux,B)$ of $\A$ such that the Cartan counterpart of the
  matrix $B$ is a Cartan matrix of finite type.
\end{theorem}
In the Theorem, the Cartan counterpart of a matrix $B=(b_{ij})$ in
$M_n(\Z)$ is the matrix $A=(a_{ij})$ with
\begin{equation}
a_{ij}=\left\{
    \begin{array}{ll}
      2 &\textup{ if } i=j,\\
      -\mid b_{ij}\mid &\textup{ if }i\not=j.
    \end{array}\right.
\end{equation}
Actually, the Theorem of Fomin and Zelevinsky is more precise. The
correspondence $\A\mapsto \Delta$ of the Theorem provides a bijection
from the set of finite reduced cluster algebras into the set of Dynkin
diagrams of finite type. Hence, to a Dynkin diagram $\Delta$ of type A
to G, we can associate a unique algebra $\ad$; this is the reduced
cluster algebra of the corresponding type.\par
Let $\cld$ be the
set of cluster variables of $\ad$. Then, by the Fomin and Zelevinsky's theorem, there exists a
  bijection $\beta$ from the set $\cld$ to $\phap$ that sends $u_i$ to
  $-\alpha_i$. This bijection will be made more precise in \ref{clustercat}.
\subsection{}
We present in this section a recollection on quiver representations.
We fix a field $k$ which can be either the finite field $\fq$ or the
field $\C$ of complex numbers. From now on, let $\Delta$ be a Dynkin
diagram of simply laced finite type and $Q$ be a quiver with
underlying graph $\Delta$. We index by $I=\{1,\ldots,n\}$ the set of
its vertices. We consider the category $\mkq$ of finite dimensional
$k$-representations of the quiver $Q$: the objects of $\mkq$ are
tuples of finite dimensional vector spaces $(M_i)_{i\in I}$ together
with tuples of linear maps $(M_\alpha:\;M_i\rightarrow
M_j)_{\alpha:\;i\rightarrow j}$, a morphism between the objects
$(M_*)$ and $(N_*)$ is a $I$-family of linear maps $M_i\rightarrow
N_i$ such that for any $\alpha:$ $i\rightarrow j$ the diagram
 \[
\xymatrix{M_i \ar[r]^{M_\alpha}\ar[d] &M_j\ar[d]  \\
N_i \ar[r]^{N_\alpha} &N_j }
\] 
commutes. We denote by $\imkq$, resp.  $\indkq$, the set of
isoclasses, resp. indecomposable modules, of $\mkq$. We also denote by
$P_i$, resp. $S_i$, resp. $I_i$, $1\leq i\leq n$, the projective,
resp. simple, resp. injective modules of the category. Recall that, by
the Theorem of Gabriel, the set $\indkq$ does not depend on the field
$k$ and is in natural bijection with the set of positive roots
$\Phi_+$.  This bijection sends the simple modules $S_i$ to the simple
roots $\alpha_i$. Hence, the Grothendieck group of the category $\mkq$
is naturally identified with the lattice $\qc$. We define the map
$\dimv$ : $\mkq\rightarrow\qc$ which associate to a module $M$ its
class $[M]$ in the Grothendieck group.
\par
We define for $M$, $N$ in $\mkq$:
\begin{equation}
[M,N]=\dim\hom_{kQ}(M,N),\,\,[M,N]^1=\dim\ext_{kQ}(M,N),
\end{equation}
\begin{equation}
<M,N>=[M,N]-[M,N]^1.
\end{equation}
It is known that $<\,,\,>$ can be defined on the Grothendieck group
$\qc$, this is the Euler form.\par
Fix $\ue$ in $\qc$ and let $M$ be in $\mkq$. We define the
$\ue$-Grassmannian $\gremk$ of the module $M$ on k:
\begin{equation}
\gremk:=\{N,\, N\in\mkq,\, N\subset M,\, \dimv(N)=\ue\}.
\end{equation}
When the fixed field is clear, we will omit the index $k$ in the
notation. Note that $\grem$ can be realized as a closed subvariety of
the classical Grassmannian $\Gr_e(M)$, where $e=\sum_i e_i$. Hence,
the variety $\grem$ is projective. We define in an obvious way the
variety $\grem$, for $M$ in $\imkq$, and we set
$\Gr(M):=\bigcup_{\ue}\grem$. So, in the sequel the Grassmannian will
always mean ``Grassmannian of submodules''.

\subsection{}\label{clustercat}
The cluster category $\cc=\ccd$ has been introduced in \cite{BMRRT}
and \cite{CCS}. As proved in these articles, the cluster category
$\ccd$ is strongly linked with the finite cluster algebra $\ad$. As a
first example of this relation, the set of indecomposable objects of
$\ccd$ is in bijection with $\phap$, so it is in bijection with the
cluster variables of $\ad$.\par
Let $Q$ be as above. In the sequel, we denote by $B_Q$ be the
antisymmetric matrix in $M_n(\Z)$ such that $b_{ij}=1$ if
$i\rightarrow j$ in $Q$ and $0$ if $i$ and $j$ are not connected. Let
${\mathcal D}\mkq$ be the derived category of $\mkq$. Note that as
$\mkq$ is hereditary, the indecomposable objects of ${\mathcal D}\mkq$
are the shifts of $\indkq$. We define the functor $F$ of ${\mathcal
  D}\mkq$ by $F$ : $M\mapsto\tau^{-1}SM$, where $S$ is the shift and
$\tau$ is the Auslander-Reiten translation.  The category $\ccq$ is
the category of orbits of ${\mathcal D}\mkq$ by $F$, see \cite{BMRRT}.
It is a triangulated category \cite{keller}, it is also a
Krull-Schmidt category, but it is not abelian in general. A nice
property of this category is that the bifunctor $\ext$ is symmetric.
This category does not depend on the orientation of $Q$, but only on
$\Delta$. It will be denoted by $\ccd$ or just $\cc$. The set $\indkq$
embeds naturally in the set $\Ind\,\cc$ of indecomposable objects of
$\cc$. Moreover, $\Ind\,\cc=\indkq\cup\{SP_i, 1\leq i\leq n\}$. Now,
via the theorem of Gabriel, the denominator theorem, \cite{FZ2},
provides a bijection between Cl$_\Delta$ and $\Ind\,\cc$.
\begin{theorem}
  Let $Q_{alt}$ be an alternating quiver with underlying diagram
  $\Delta$. Consider the seed $(\underline u=\{u_1,\ldots,u_n\},
  B_{Q_{alt}})$ and a cluster variable $u$, $u\not=u_i$, $1\leq i\leq
  n$. Then, there exists a positive root $\alpha=\sum_in_i\alpha_i$
  such that the denominator of $u$ as an irreducible fraction in
  $\underline u$ is $\prod_iu_i^{n_i}$. The assignment
  $u\mapsto\alpha$, $u\not=u_i$, $u_i\mapsto -\alpha_i$, provides a
  bijection from Cl$_\Delta$ to $\phap$.
\end{theorem}

An important
Theorem of \cite{BMRRT} asserts that the ext-configurations of
$\Ind\,\cc$, \textit{i.e.} maximal subsets of $\Ind\,\cc$ with trivial
pairwise extension, correspond to clusters via the bijections
Ind$(\ccd)\simeq\phap\simeq\cld$. Note also that the
ext-configurations of $\indkq$ are ext configurations of $\Ind\,\cc$
which belong to $\indkq$.

\subsection{}\label{wheels}
In the previous section, we have seen a correspondence between
indecomposable objects of $\cc$ and cluster variables. We will see in
this section properties of the Auslander-Reiten translation $\tau$ in
this correspondence. First of all, let us recall some basic facts on
the Auslander-Reiten theory, see \cite{ARS}.\par
Let $\gq$ be the Auslander-Reiten quiver of $\mkq$. Recall that its
set of vertices is $\indkq$ and the arrows are given by irreducible
morphisms of the category. The AR-quiver $\gcc$ of $\cc$ is defined in
the same way.\par
Let $M$ be a non projective module in $\indkq$ and let $M'=\tau
M\in\indkq$ be its AR-translated. We consider the direct sum $B$ of
indecomposable modules $B_j$ such that $M'\rightarrow B_j$ in $\gq$.
Then, $B$ is also the direct sum of indecomposable modules $B_j$ such
that $B_j\rightarrow M$ in $\gq$ and we have the following exact
sequence of modules:
\begin{equation}\label{asseq}
\xymatrix{
0\ar[r] &\tau M\ar[r] &B
\ar[r]^{\sigma} &M\ar[r] &0}.
\end{equation}
                              
Moreover, this exact sequence is almost split in the following sense:
each morphism $N\rightarrow M$ which is not a split epimorphism
factors through $\sigma$. \par
The AR-quivers of $\mkq$ are well known and can be explicitly
described, see \cite{gabriel}.\par
The AR-quiver $\gcc$ is a slight extension of $\gq$. Indeed, see
\cite{BMRRT}, each exact sequence as in (\ref{asseq}) gives rise in
the triangulated category $\cc$ to a triangle
\begin{equation}    
\tau M\rightarrow B\rightarrow M\rightarrow S\tau M,
\end{equation}
where the first two morphisms are composed with irreducible morphisms.
In other words, the embedding $\indkq\subset\Ind\,\cc$ provides an
embedding $\gq\subset\gcc$ as a full subquiver. In order to describe
$\gcc$, it is sufficient to note that:\par
There exists an arrow $SP_i\rightarrow M$ in $\cc$ if and only if
\begin{equation}\label{dis1}
  M=\left\{\begin{array}{ll}
      SP_k & k\rightarrow i\textup{ in } Q,\\
      P_j & i\rightarrow j\textup{ in } Q,
    \end{array}\right.
\end{equation}
and there exists an arrow $M\rightarrow SP_i$ in $\cc$ if and only if
\begin{equation}\label{dis2}
  M=\left\{\begin{array}{ll}
      P_k & k\rightarrow i\textup{ in } Q,\\
      SP_j & i\rightarrow j\textup{ in } Q.
    \end{array}\right.
\end{equation}

In the following proposition, we denote by $x_M$ the cluster variable
corresponding to the indecomposable object $M$ of $\Ind\,\cc$.
    
\begin{proposition}
  Let $\cc$ be the cluster category of type A, D or E.
  \begin{itemize}
  \item[(i)] Let $M$ be an indecomposable object of $\cc$. Then $x_M$
    and $x_{\tau M}$ form an exchange pair.
  \item[(ii)] With the notation above, we have 
    \begin{equation}
      x_{\tau M}x_M=\prod_j x_{B_j}+1.
    \end{equation}
  \end{itemize}
\end{proposition}
\begin{proof}
  (i) is a direct application of \cite[Proposition 7.6]{BMRRT}.\par
  Fix $M$ in $\Ind\,\cc$. Then, there exists a quiver $Q$ and a sink
  $i$ of $Q$ such that in the equivalence $\cc\simeq {\mathcal
    D}\mkq\slash F$, the object $M$ is identified with the simple
  projective $P_i=S_i$. The algebra $\A(\Delta)$ is by construction
  isomorphic to $\A(B_Q)$. By the discussion above, see (\ref{dis1})
  and (\ref{dis2}), it is enough to prove that
  \begin{equation}
    x_{S P_i}x_{S_i}=\prod_{i\rightarrow j} x_{S P_j}+1. 
  \end{equation}
  Let us prove the equality. Set $B_Q=(b_{ij})$. The exchange relation
  gives
  \begin{equation}
    x_{S P_i}x_{S_i}=\prod_{b_{ji}=1} x_{S P_j}+  \prod_{b_{ji}=-1}
    x_{S P_j}.
  \end{equation}
  As $i$ is a sink of $Q$, the second term is one. Moreover,
  $b_{ji}=1$ if and only if $j\rightarrow i$ in $Q$. So, we have the
  claimed equality.
  \end{proof}
  Remark that (ii) is a particular case of \cite[Conjecture 9.3]{BMRRT}.
\section{The main Theorem}
\subsection{}
In this section, we give realizations of any finite reduced cluster
algebra from the category $\cc$; for each quiver $Q$ with
underlying Dynkin diagram $\Delta$ of type A-D-E, we realize $\ad$
from $\ccq$. Actually, we recover the algebra $\A(B_Q)$ from the
category $\ccq$. In the sequel, we fix a quiver $Q$ of type A-D-E.\par
Recall that $\F=\Q(u_i, 1\leq i\leq n)$. For each $M$ in $\imkq$ with
dimension vector $\dimv(M)=\um=\sum_i m_i \alpha_i$, set
\begin{equation}  \label{explicit}
  X_M=\sum_{\ue}\chi(\grem)\prod_i u_i^{-<\ue,\alpha_i>-<\alpha_i,\um-\ue>},
\end{equation}
where $\chi$ is the Euler-Poincaré characteristic of the complex
Grassmannian. Remark that the sum is finite since the dimension
vectors $\ue=\sum_i e_i\alpha_i$ which occur in the sum verify $0\leq
e_i\leq m_i$. We now illustrate with examples.
\begin{example}\label{A3}
  Suppose that $Q$ is the following alternated orientation for A$_3$:
  \begin{equation}
    \xymatrix{
      1 \ar[r] & 2 & 3 \ar[l]
      }.
  \end{equation}
  Then, the indecomposable modules of $\mkq$ are $S_1$, $S_2$, $S_3$,
  $P_1$, $P_3$, $I_2$, where $[I_2]=[S_1]+[S_2]+[S_3]$. The AR-quiver
  $\gcc$ has the following shape:
  \begin{equation}
    \xymatrix{
       & S P_3 \ar[rd] & & P_3 \ar[rd] & & S_1 \ar[dr] & \\
       S P_2 \ar[ru] \ar[rd] & & S_2 \ar[ru] \ar[rd] & & I_2 \ar[ru]
       \ar[rd] & & S P_2\\
       & S P_1 \ar[ru] & & P_1 \ar[ru] & & S_3 \ar[ur] &.\\
     }
  \end{equation}
  We compute explicitly the $X_M$ using formula (\ref{explicit}). In
  the following sums, the terms are ordered by $\sum e_i$.
  The submodules of $S_2$ are $0$ and $S_2$.
  \begin{equation}
    X_{S_2}=\frac{u_1u_3}{u_2}+\frac{1}{u_2}=\frac{u_1u_3+1}{u_2},
  \end{equation}
  The submodules of $P_3$ are $0$, $S_2$ and $P_3$ and the submodules
  of $P_1$ are $0$, $S_2$ and $P_1$.
  \begin{equation}
    X_{P_3}=\frac{u_1}{u_2}+\frac{1}{u_2u_3}+\frac{1}{u_3}=\frac{1+u_2+u_1u_3}{u_2u_3},\,
    X_{P_1}=\frac{u_3}{u_2}+\frac{1}{u_2u_1}+\frac{1}{u_1}=\frac{1+u_2+u_1u_3}{u_2u_1},
  \end{equation}
  The submodules of $I_2$ are $0$, $P_1$, $P_3$, $S_2$ and $I_2$.
  \begin{equation}
    X_{I_2}=\frac{1}{u_2}+\frac{1}{u_1u_3}+\frac{1}{u_1u_3}+\frac{1}{u_1u_2u_3}+\frac{u_2}{u_1u_3}=
    \frac{1+2u_2+u_2+u_1u_3}{u_1u_2u_3},
  \end{equation}
  The submodules of $S_1$ are $0$ and $S_1$ ; the submodules of $S_3$
  are $0$ and $S_3$.
  \begin{equation}
    X_{S_1}=\frac{1}{u_1}+\frac{u_2}{u_1}=\frac{1+u_2}{u_1},\;
    X_{S_3}=\frac{1}{u_3}+\frac{u_2}{u_3}=\frac{1+u_2}{u_3}.
  \end{equation}
\end{example}
\begin{example}\label{An}
  If $Q$ is a quiver of type A$_n$, and if $M$ is an indecomposable
  module of $\mkq$, then $\chi(\grem)=0$ or $1$. More precisely, the
  indecomposable $kQ$-modules correspond to connected full subquivers
  of $Q$. Let $Q_M$ be the quiver corresponding to $M$ and let $V_M$
  be the set of its vertices. Then, the submodules $N$ of $M$
  correspond to subsets $V_N$ of $V_M$ such that the following
  property holds: $ i\in V_N$ and $i\rightarrow j\Rightarrow j\in
  V_N$.\par
  Now, if $N$ is a submodule of $M$ with dimension vector $\un$, then
  $\Gr_{\un}(M)$ has only one point and $\chi(\Gr_{\un}(M))=1$.
  
  In the particular case where $Q$ is the equioriented quiver of type
  A$_n$, the property above implies that, for each indecomposable
  module $M$, the number of terms in the decomposition of $X_M$ is
  dim$M+1$.
\end{example}
\begin{example}
  We consider the following quiver $Q$ of type D$_4$:
  \begin{equation}
    \xymatrix{
      & 3 \ar[d] & \\
      1 \ar[r] & 2 & \ar[l] 4.
    }
  \end{equation}
  Let $M$ be the indecomposable module with maximal dimension,
  \textit{i.e.} $[M]=[S_1]+[S_3]+[S_4]+2[S_2]$. Then, we have
  $\Gr_{\alpha_2} M={\mathbb P}^1$ and so $\chi(\Gr_{\alpha_2} M)=2$.
  \par
  The module $M$ has 13 submodules: $0$, $S_2$, $2S_2$, $P_1$, $P_3$,
  $P_4$, $P_1+S_2$, $P_3+S_2$, $P_4+S_2$, $P_1+P_3$, $P_1+S_2$,
  $P_1+P_4$, $P_3+P_4$, $M$. But, $S_2$ has "multiplicity" 2. That
  gives
  \begin{equation}
    X_M=\frac{ (1+ u_2)^3 + 2 u_1 u_3 u_4 + 3 u_1 u_2 u_3 u_4  
        + u_1^2  u_3^2  u_4^2 }{u_1 u_2^2 u_3 u_4}.
  \end{equation}
\end{example}

Let $E_Q$ be the $\Q[u_i,\,1\leq i\leq n]$-submodule of $\F$ generated
by $X_M$, $M\in\imkq$, then:
\begin{theorem}\label{main}
  For each quiver $Q$ of type A-D-E, $E_Q$ is a subalgebra of $\F$. It
  identifies with the subalgebra $\A(B_Q)=\ad$ of $\F$. Up to this
  identification, the set of cluster variables of $\A(B_Q)$ is given
  by $\{u_i,\,1\leq i\leq n\}\cup\{ X_M,\,M\in\indkq\}$.
\end{theorem}
Note that in particular, up to isomorphism, the algebra $E_Q$ does not
depend on $Q$ but only on $\Delta$.



The subsections below are devoted to the proof of this Theorem.

\subsection{}
In order to calculate the Euler-Poincaré characteristic of
Grassmannians, we will use the following classical Lemma, see \cite{rei}, which is an
application of Grothendieck-Lefschetz's fixed point formula for the
Frobenius in {\'e}tale cohomology.
\begin{lemma}\label{deligne}
  Let $X$ be a variety defined over some ring of algebraic integers.
  We denote by $X_{\C}$ (resp. $X_{\fq}$) the set of $\C$-points (
  resp. $\fq$-points) of $X$. Suppose that there exists a polynomial
  $P$ with integral coefficients such that $\mid X_{\fq}\mid=P(q)$ for
  infinitely many prime powers $q$. Then, the Euler-Poincaré
  characteristic (with compact support) of $X_{\C}$ is given by
  $\chi(X_{\C})=P(1)$.
\end{lemma}

Note first that the Grassmannians discussed above are defined over
$\Z$. Indeed, Grassmannians of $k$-subspaces are defined by base
change from a $\Z$-scheme, \cite[3.13]{DG}. Moreover, as the
$kQ$-modules are defined on $\Z$, it is easily seen that a subspace is
a submodule if and only if it verifies $\Z$-linear conditions on the
Pl\"ucker coordinates of the subspace. Note also that the cardinality
of Grassmannians of submodules is given by sums of Hall polynomials.
Indeed, for a fixed module $M$ and a fixed dimension vector of a
submodule $N$, there is only a finite number of possibilities for the
isomorphism classes of $N$ and of the quotient $M/N$, as their
dimension vectors are fixed and the quiver is of finite type. Hence
the cardinality of a Grassmannian of submodules of fixed dimension
vector is a finite sum of the cardinalities of sets of triples
$(N,M,M/N)$ where the isomorphism classes are fixed. These
cardinalities are known to be polynomials in $q$, called Hall
polynomials, see \cite{ringel}.

\subsection{}\label{subalgebra}
We prove here that $E_Q$ is a subalgebra of $\F$. Actually, we will
prove the following:
\begin{proposition}\label{split}
  Fix $\ug$ in $\qc$. For all $M$, $N$ in $\imkq$, we have
  \begin{equation}
    \chi(\Gr_{\ug}(M\oplus N))
    =\sum_{\ue+\uf=\ug}\chi(\Gr_{\ue}(M))\chi(\Gr_{\uf}(N)).
  \end{equation}
\end{proposition}
By the bilinearity of the Euler form, this Proposition implies
\begin{corollary}\label{algebra}
  For all $M$, $N$ in $\imkq$, we have $X_M X_N=X_{M\oplus N}$. Hence,
  $E_Q$ is a subalgebra of $\F$. It is the $\Q$-subalgebra generated
  by $\{u_i,\,1\leq i\leq n\}\cup\{ X_M,\,M\in\indkq\}$.
\end{corollary}
By Lemma \ref{deligne}, Proposition \ref{split} can be obtained by
counting points on $\fq$ varieties. Set $k=\fq$. Fix two $kQ$-modules
$M$ and $N$. Let $\pi$ : $M\oplus N\rightarrow N$ be the projection on
the second factor.

Fix a submodule $A$ of $M$ and a submodule $B$ of $N$. Let us
introduce
\begin{equation}
  \Gr_{A,B}(M\oplus N):=\{L\in\Gr(M\oplus N),\,L\cap M=A\hbox{ and }\pi(L)=B\}.
\end{equation}

\begin{lemma}
  Fix a submodule $A$ of $M$ and a submodule $B$ of $N$. Let $\pi_A$ : $M\rightarrow A$ be the 
  canonical projection. There
  exists a bijection:
  \begin{equation}
    \hom_{kQ}(B,M/A)\rightarrow \Gr_{A,B}(M\oplus N), 
  \end{equation}
  which maps the morphism f to $L_f=\{m+b,\,m\in M,\,b\in B,\,\pi_A(m)=f(b)\}$.
\end{lemma}
\begin{proof}
  The space $L_f$ is a submodule of $M\oplus N$ and $\pi(L_f)=B$ by
  construction. Moreover, $L_f\cap M=$Ker$(\pi_A)=A$, hence,
  $L_f\in\Gr_{A,B}(M\oplus N)$ and the map is well defined.\par
  To show that the correspondence $f\mapsto L_f$ is bijective, define
  the opposite direction map as follows: let $L\in \Gr_{A,B}(M\oplus
  N)$. For $b$ in $B$, define $f_L(b)=\pi_A(m)$ for any $m\in M$ such
  that $m+b\in L$. Since $L\cap M=A$, this map is well defined.
  Clearly, both left and right compositions with $f\mapsto L_f$ are
  identity maps, as needed.
 
\end{proof} 
{\it Proof of the Proposition.}
For $\ug$ in $\qc$, consider the map
\begin{equation}
  \zeta_{\ug}\,:\,\Gr_{\ug}(M\oplus N)\rightarrow 
  \coprod_{\ue+\uf=\ug}
  \Gr_{\ue}(M)\times\Gr_{\uf}(N),\,L\mapsto (L\cap M, \pi(L)).
\end{equation}
This map is clearly surjective: $\zeta_{\ug}(A\oplus B)=(A,B)$.
Moreover, the Lemma above proves that $\zeta_{\ug}^{-1}(A,B)$ has
$q^{[B,M/A]}$ elements. Now, the Proposition is a direct consequence
of Lemma \ref{deligne}.

\subsection{}\label{generators}
Now, we need to understand the natural set of generators of the
algebra $E_Q$, $\{u_i,\,1\leq i\leq n\}\cup\{ X_M,\,M\in\indkq\}$ by
Corollary \ref{algebra}. We want to prove that it is precisely the set
of cluster variables of $\A(B_Q)$. By construction the variables
$u_i$, $1\leq i\leq n$, are cluster variables. Now, for each
indecomposable module $M$ in $\mkq$, let $\nu(M)$ be the smallest
integer such that $\tau^{\nu(M)}M=0$ in $\mkq$. We want to prove by
induction on $\nu(M)$ that $X_M$ is a cluster variable of $\A(B_Q)$.
The case $\nu(M)=1$ corresponds to the projective case. By Section
\ref{wheels}, in this case, we have to prove:
\begin{lemma}
  For all $i$, $1\leq i\leq n$, we have
  \begin{equation}
    u_iX_{P_i}=\prod_{i\rightarrow j}X_{P_j}\prod_{k\rightarrow i}u_k+1. 
  \end{equation}
\end{lemma}
\begin{proof}
  Set $\ud_i:=\dimv(P_i)$. It is known that the radical $\Rad P_i$
  verifies the following:
  \begin{itemize}
  \item[(i)] $P_i/\Rad P_i=S_i$, 
  \item[(ii)] $M\subset P_i\Leftrightarrow M\subset\Rad P_i$ or $M=P_i$,
  \item[(iii)] $\Rad P_i=\oplus_{i\rightarrow j}P_j$.
  \end{itemize}
  By (i), we have
  \begin{equation}
    X_{\Rad P_i}=\sum_{\ue}\chi(\Gr_{\ue}\Rad P_i)\prod_l u_l^{-<\ue,\alpha_l>-<\alpha_l,\ud_i-\alpha_i-\ue>}.
  \end{equation}
  Using the fact that the Euler form satisfies
  \begin{equation}
    < \alpha_k, \alpha_i >=
    \begin{cases}
      1 & \text{ if }k=i,\\
      -1 & \text{ if }k\rightarrow i,\\
    0 & \text{ else,}
    \end{cases}
  \end{equation}
  one gets
  \begin{equation}\label{firsteq}
    X_{\Rad P_i}=\sum_{\ue}\chi(\Gr_{\ue}\Rad P_i)(\prod_l u_l^{-<\ue,\alpha_l>-<\alpha_l,\ud_i-\ue>})(\prod_{k\rightarrow i}u_k^{-1})u_i.
  \end{equation}
  By (ii), we have
  \begin{equation}
    X_{P_i}=\sum_{\ue}\chi(\Gr_{\ue}\Rad P_i)\prod_l u_l^{-<\ue,\alpha_l>-<\alpha_l,\ud_i-\ue>}+u_i^{-1}.
  \end{equation}
  Comparing with (\ref{firsteq}) gives
  \begin{equation}
    X_{P_i}=X_{\Rad P_i}(\prod_{k\rightarrow i} u_k)u_i^{-1}+u_i^{-1}.
  \end{equation}
  The Lemma is now a consequence of (iii) and Corollary \ref{algebra}.
\end{proof}

\subsection{}
We prove here the induction discussed in Section \ref{generators}.
What we need to prove is that for all non projective indecomposable
$kQ$-module $N$, if $X_{\tau N}$ is a cluster variable of
$\A(\Delta)$, then $X_N$ is also a cluster variable. By Section
\ref{wheels}, what we have to prove is
\begin{proposition}\label{relations}
  Suppose that $M$, $N$ are indecomposable modules and 
  \begin{equation}\label{arseq}
    \xymatrix{0\ar[r] &M\ar[r]^{\iota} &B \ar[r]^{\pi} &N\ar[r] &0}
  \end{equation}
  is an almost split exact sequence, then $X_{M\oplus
    N}=X_M X_N=X_B+1$.
\end{proposition}
Remark that $M\oplus N$ and $B$ are the middle terms $Y$ of
respectively a split sequence and an almost split sequence
$0\rightarrow M\rightarrow Y\rightarrow N\rightarrow 0$. The reader
may view the Formula above as a ``difference'' between split and
almost split. The proof of the Proposition is an adaptation of the
proof in Section \ref{subalgebra} in the almost split case.
\begin{proof}
  Set $\um=\dimv M$, $\un=\dimv N$. Recall that $M=\tau N$. We have
  \begin{equation}
    X_{M\oplus N}=
    \sum_{\ue}\chi\left(\Gr_{\ue}(M\oplus N)\right)
    \prod_i u_i^{-<\ue, \alpha_i>-<\alpha_i,\um+\un-\ue>}. 
  \end{equation}
  By Lemma \ref{deligne}, $X_{M\oplus N}$ can be seen as a polynomial
  of $\Z[u_i^{\pm 1}][q]$ evaluated at $q=1$:
  \begin{equation}
    X_{M\oplus N}=(\sum_L\prod_i u_i^{-<\dimv L, \alpha_i>-<\alpha_i,\um+\un-\dimv L>})\mid_{q=1}, 
  \end{equation}
  where $L$ runs over the set of submodules of $M\oplus N$ and
  $k=\fq$. In this Formula the term corresponding to the submodule
  $L=0\oplus N$ in $M\oplus N$ is
  \begin{equation}
    \prod_i u_i^{-<\un,\alpha_i>-<\alpha_i,\um>}=1,
  \end{equation}
  by the Serre duality formula.
  \par
  As in the proof of Proposition \ref{split}, our Proposition follows
  from Lemma below.
\end{proof}
\begin{lemma}\label{merciK}
  Consider the map
  \begin{equation}
    \zeta_{\ug}\,:\,\Gr_{\ug}(B)\rightarrow 
    \coprod_{\ue+\uf=\ug}\Gr_{\ue}(M)\times\Gr_{\uf}(N),
    \,L\mapsto (\iota^{-1}(L), \pi(L)),
  \end{equation}
  where $\iota$ and $\pi$ are the morphisms in \ref{arseq}.
  The fiber of a point $(A,C)$ is empty if $(A,C)=(0,N)$, and is an
  affine space of dimension $[C,M/A]$ if not.
\end{lemma}
\begin{proof}
  Let us prove the case $(A,C)=(0,N)$. Suppose that $L\subset B$,
  $\pi(L)=N$ and $\iota^{-1}(L)=0$. Then, $\pi$ provides an
  isomorphism $L\simeq N$. This implies that the map $\pi$ splits but
  this is impossible since $\pi$ is the surjection of an almost split
  sequence. \par
  Suppose now $(A,C)\not=(0,N)$. If $C$ is not equal to $N$, then, the
  ``almost split'' property implies that $\pi$ has a section
  $C\rightarrow B$. We are in the split case and the proof is as in
  Section \ref{subalgebra}.\par
  It remains to prove the case where $A\not=0$ and $C=N$. Since
  $0\not=A\subset M$, we have $[N,A]^1=[A,M]\not=0$. There exists a
  non split exact sequence $0\rightarrow A\rightarrow E\rightarrow
  N\rightarrow 0$. We claim that we have following commutative
  diagram:
  \begin{equation}
    \xymatrix{{}&0\ar[d]&0\ar[d]&0\ar[d]&{}\\ 
      0\ar[r]&A\ar[d]\ar[r]&E\ar[d]^{\varphi}\ar[r]_{\mu}
      &N\ar@{=}[d]\ar[r]&0\\
      0\ar[r]&M \ar[r]& B \ar[r]^{\pi} & N\ar[r]& 0.}
  \end{equation}
  Indeed, $\varphi$ exists by the almost split property and it is a
  monomorphism by an easy diagram chasing.\par
  Hence, $E\in\zeta_{\ug}^{-1}(A,N)$ and the fiber is non empty. In
  this case, it is a well-known fact that the fiber is an affine space
  with a simple transitive action of $\hom(N,M/A)$. We give the proof
  for completion. \par
  Let $\mu$ : $E\rightarrow N$ as in the diagram above and $\pi_A:$
  $M\rightarrow M/A$ the canonical surjection. We consider the map
  $\hom(N,M/A)\rightarrow \zeta_{\ug}^{-1}(A,N)$, defined by $f\mapsto
  E_f:=\{\iota(m)+\varphi(e),\,m\in M,\, e\in E,\,\pi_A(m)=f\mu(e)\}$.
  The space $E_f$ is a submodule of $B$. We start with the following
  remark. Suppose that
  $$\iota(m')+\varphi(e')=\iota(m)+\varphi(e),\,m,\,m'\in M,\,e,\,e'\in E.$$
  Applying $\pi$ gives $e'-e\in\iota(A)$ and then
  $\pi_A(m)=\pi_A(m')$. In particular, this easily implies that $E_f$
  is in the fiber $\zeta_{\ug}^{-1}(A,N)$.\par
  
  In order to prove that the correspondence $f\mapsto E_f$ is
  bijective, define the opposite direction map as follows: let
  $D\in\zeta_{\ug}^{-1}(A,N)$, for $n$ in $N$, set
  $f_{D}(n)=\pi_A(m)$, where $\iota(m)+\varphi(e)\in D$ and
  $\mu(e)=n$. By the remark above, $f_{D}$ is a well defined element
  of $\hom(N,M/A)$.
  
  We have $E_{f_{D}}=D$. Indeed, it is enough to prove the inclusion,
  as both modules have same dimension vector $\underline g$. Let $x$
  in $E_{f_{D}}$, hence we have a decomposition
  $x=\iota(m)+\varphi(e)$, with $\pi_A(m)=f_{D}(\mu(e))$. By
  construction of $f_{D}$, we have $\pi_A(m)=f_{D}\mu e=\pi_A(m')$,
  with $\iota(m')+\varphi(e)\in D$. Hence,
  $x=\iota(m)+\varphi(e)=\iota(m')+\varphi(e)\in D$, as desired.
  
  We have $f_{E_f}=f$. Indeed, $f_{E_f}(n)=\pi_A(m)$, where
  $\iota(m)+\varphi(e)\in E_f$ and $\mu(e)=n$. By definition of $E_f$
  we have: $\iota(m)+\varphi(e)=\iota(m')+\varphi(e')$, and
  $\pi_A(m')=f\mu(e')$. So, by the remark above,
  $$f_{E_f}(n)=\pi_A(m)=\pi_A(m')=f\mu(e')=f\mu(e)=f(n).$$ This ends the proof.
 \end{proof}

\section{A conjecture in the multiplicity-free case}

\subsection{}

Let $Q$ be a Dynkin quiver. In particular the category of modules over
$Q$ is hereditary. Let $M$ be an indecomposable object in the category
of modules over the quiver $Q$. Assume that $M$ is multiplicity-free,
that is $\dim(M_i)\leq 1$ for all $i$. This implies that a submodule
$N$ of $M$ is determined by its dimension vector. In these cases, the
Grassmannian is either empty or a point.

Then Formula (\ref{explicit}) for the cluster corresponding to $Q$ and
the cluster variable corresponding to $M$ can be restated as follows.

\begin{proposition}
  One has 
  \begin{equation}
    \label{multfree}
    X_M=\frac{1}{\prod_{i\in M}u_i}
    \sum_{N \subset M} {\prod_{i \in N}\left(\prod_{i
          \to j} u_j \right)\prod_{i \in M/N}\left(\prod_{j\to i} u_j
      \right)},
  \end{equation}
  where the sum runs over submodules $N$ of $M$ and the index $j$ in
  the inner products runs over the set of vertices of $Q$.
\end{proposition}

\begin{proof}
  Indeed, Formula (\ref{explicit}) can be reformulated, using
  injective and projective resolutions, the definition of the Euler
  form and the known Euler-Poincaré characteristic of Grassmannians,
  as the following expression:
  \begin{equation}
    X_M=\sum_{N \subset M}
    \frac{[P^N_1]}{[P^N_0]}\frac{[I^{M/N}_1]}{[I^{M/N}_0]},
  \end{equation}
  where
  \begin{equation}
    \xymatrix{
      0 \ar[r] & P^N_1 \ar[r] & P^N_0\ar[r] & N \ar[r] & 0}
  \end{equation}
  and
  \begin{equation}
    \xymatrix{
      0 \ar[r] &  M/N  \ar[r]& I^{M/N}_0 \ar[r] & I^{M/N}_1 \ar[r] & 0}
  \end{equation}
  are projective and injective resolutions and the brackets mean
  replacing the direct sum of projective modules $P_i$ or injective
  modules $I_i$ by the corresponding product of variables $u_i$.

  Let us fix a submodule $N$ of $M$ and denote by $\source(N)$ the set
  of sources of the quiver underlying $N$.

  Using the hypotheses that the module $N$ is multiplicity-free and
  that the quiver $Q$ is Dynkin of finite type, hence a tree, one can
  describe completely the minimal projective resolution of $N$. The
  first step $P^N_0$ is the direct sum of all projective modules $P_j$
  for $j \in \source(N)$. By the hereditary property, the dimension
  vector of the second step $P^N_1$ is known. Its support is made of
  some isolated elements inside $N$ and some branches starting just
  outside $N$.  Each relative sink $j$ in $N\setminus \source(N)$
  contributes to $P^N_1$ by the direct sum of $\mathsf{N}(j)-1$ copies
  of $P_j$ where $\mathsf{N}(j)$ is the number of arrows inside $N$
  with target $j$.  Each of the branches corresponds to a projective
  $P_j$ for some $j$ outside $N$ with an arrow $i \to j$ for some $i
  \in N$.

  From this description of the projective resolution of $N$, one gets that
  \begin{equation}
    \frac{[P^N_1]}{[P^N_0]}=
    \frac{\prod_{i\in N}\left(\prod_{j \not\in
          N,i\to j}u_j \right) 
      \prod_{j\in N \setminus\source(N)} u_j^{\mathsf{N}(j)-1}}
    {\prod_{j \in \source(N)}u_j}.
  \end{equation}
  This becomes
  \begin{equation}
    \prod_{i\in N}\left(\prod_{j \not\in
        N,i\to j}u_j \right) 
    \prod_{j\in N} u_j^{\mathsf{N}(j)-1}.
  \end{equation}
  
  Then it follows that
  \begin{equation}
    \frac{[P^N_1]}{[P^N_0]}=
    \frac{\prod_{i \in N}\left( \prod_{i \to j}u_j\right)}{\prod_{j \in N}u_j},
  \end{equation}
  where the index $j$ in the numerator product runs over the set of
  vertices of $Q$. A similar argument for the injective resolution of
  $M/N$ completes the proof.
\end{proof}
\subsection{} 

Let us now consider a quiver $Q$ of finite cluster type, as introduced
in \cite{FZ2} and studied in \cite{CCS}. This can be one of the
quivers of Dynkin type considered before, but many other quivers arise
in the mutation process starting from a Dynkin quiver. These quivers
can be defined from the matrix $B$ of a seed $(\ux,B)$ in a
simply-laced finite reduced cluster algebra by the rule that there is
an arrow from $i$ to $j$ if and only if one has $b_{ij}=1$. This rule
was already used for the Dynkin quivers in Section \ref{clustercat}.

Then it is expected in general and known in type A \cite{CCS} that
there is a correspondence between cluster variables (other than the
initial ones) for the seed associated to $Q$ and indecomposables of
the category of modules over $Q$ with some relations. In this
bijection, the denominators of the cluster variables should be
described by the dimension vectors of the indecomposables.

Although the precise relations are not known outside of type A, a
conjecture for them has been made in \cite{CCS} and it is usually a
simple task to check in any particular case that the proposed
relations have the expected properties.

\medskip

Assuming now that the proposed relations are correct or that the
correct relations are known, let us propose a formula for the cluster
variables associated to multiplicity-free indecomposables. Let $M$ be
a multiplicity-free indecomposable object in the category of modules
over the quiver $Q$ with relations. Abusing notation, we will denote a
submodule $N$ of $M$ and its support by the same letter.

Let $\edge_M$ be the set of arrows of $Q$ between vertices of $M$ such
that the associated morphism in $M$ is zero. For $e$ in $\edge_M$, let
$s(e)$ and $t(e)$ be the source and target of $e$. We will display
later an example of cluster quiver and a module over it in type D$_4$
where $\edge_M$ is not empty.

\begin{conjecture}
  \label{mfreeconjecture}
  The cluster variable $X_M$ has the following expression:
  \begin{equation}
    \label{conj_form}
    \frac{1}{\prod_{i\in M}u_i}
    \sum_{N \subset M} \frac{\prod_{i \in N}\left(\prod_{i
          \to j} u_j \right)\prod_{i \in M/N}\left(\prod_{j\to i} u_j
      \right)}
    {\prod_{i\not\in M, i\to M, M\to i}u_i \prod_{e \in \edge_M}u_{s(e)}u_{t(e)}},
  \end{equation}
  where the sum runs over submodules $N$ of $M$. Here $i \to M$ means
  that there exists $k$ in $M$ and an arrow $i \to k$. The meaning of
  $M \to i$ is similar. The index $j$ in the inner products runs over
  the set of vertices of $Q$.
\end{conjecture}

One can recognize in the left factor of this Formula the expression
for what should be the denominator of the cluster variable. Therefore
the remaining part should be a formula for the numerator.

In the case of a Dynkin quiver $Q$, as $M$ is assumed indecomposable,
there can not be any vanishing edge in the support of $M$, hence
$\edge_M$ is empty. As $Q$ is a tree, there is no vertex outside $M$
with arrows in $Q$ to $M$ and from $M$. Hence Formula (\ref{multfree})
is a special case of Conjecture \ref{mfreeconjecture} and this
conjecture holds in the case of Dynkin quivers.

\subsection{}
Let us give two simple examples of Formula (\ref{conj_form}) for
non-Dynkin cluster quivers. Let us consider first the following
cluster quiver $Q$ of type A$_3$:
\begin{equation}
  \xymatrix{
    1 \ar[rd] &  \\
    2 \ar[u] & 3, \ar[l]
}
\end{equation}
with the relations $f_{1,2}\,f_{2,3}=f_{2,3}\,f_{3,1}
=f_{3,1}\,f_{1,2}=0$ according to \cite{CCS}. Let $M$ be the
indecomposable multiplicity-free module with support $\{1,2\}$. It has
$3$ submodules: $0,M$ and a submodule with support $\{1\}$. Then
Formula (\ref{conj_form}) gives
\begin{equation}
  X_M=\left(\frac{1}{u_1 u_2}\right)\frac{(1)(u_3 u_1)+(u_3)(u_3)+(u_3
    u_2)(1)}{u_3}
=\frac{u_1+u_2+u_3}{u_1 u_2},
\end{equation}
which is the correct expression.

Let us consider now the following cluster quiver $Q$ of type D$_4$:
\begin{equation}
  \xymatrix{
    1 \ar[r] \ar[d] & 2 \ar[d] \\
    3 \ar[r] & 4, \ar[ul]
}
\end{equation}
with relations $f_{1,2}\, f_{2,4}=f_{1,3}\, f_{3,4}$ and $f_{2,4}
\,f_{4,1}=f_{4,1}\, f_{1,2}=f_{4,1}\, f_{1,3}=f_{1,3}\, f_{3,4}=0$, as
conjectured in \cite{CCS}. One can easily compute the indecomposable
objects and the Auslander-Reiten quiver for the category of modules on
this quiver with these relations. One can then check that the
dimension vectors of the indecomposable modules correspond to the
denominators of the cluster variables when expressed in this cluster.
This proves in this particular case some conjectures made in
\cite{CCS}. Let $M$ be the following indecomposable module:
\begin{equation}
  \xymatrix{
    k \ar[r]^{\id} \ar[d]^{\id} & k \ar[d]^{\id} \\
    k \ar[r]^{\id} & k. \ar[ul]_0
}
\end{equation}
Note that the set $\edge_M$ contains the diagonal arrow.

The module $M$ has $6$ submodules: $0,M$ and submodules with support
$\{4\},\{2,4\},\{3,4\}$ and $\{2,3,4\}$. Then using Formula
(\ref{conj_form}), one gets that $X_M$ is equal to
\begin{equation*}
  \frac{(1)(u_1^2 u_2 u_3 u_4)+(u_1)(u_1^2 u_4)+(1+1)(u_1 u_4)(u_1 u_4)
    +(u_1 u_4^2)(u_4)+(u_1 u_2 u_3 u_4^2)(1)}{(u_1 u_2 u_3 u_4)(u_1 u_4)}
,
\end{equation*}
which simplifies to the correct expression:
\begin{equation}
  X_M
  =\frac{(u_1+u_4)^2+u_2 u_3 (u_1+u_4)}{u_1 u_2 u_3 u_4}.
\end{equation}
\section{Coxeter-Conway friezes}
We give here an interpretation of Coxeter-Conway friezes
\cite{CoxCon}, which follows directly from the main theorem.
\subsection{}
Following Conway and Coxeter, we construct a frieze from a
triangulation of the $(n+3)$-gon. The construction is the
following.\par
We consider a triangulation $T$ of the $(n+3)$-gon, \textit{i.e.} a
maximal set of non crossing diagonals of the polygon. Note that each
maximal set has exactly $n$ diagonals. To each vertex $k$,
$k\in\Z/(n+3)\Z$, of the polygon, let $d_k$ be the number of diagonals
of $T$ containing the vertex $k$. Set $m_k=d_k+1$. We construct a
frieze filled with numbers in the following way.\par
We place a bottom row R$_0$ filled with 1. Then, we place above it (in
a shifted way, see Example \ref{frieze} below) a row filled with
$m_k$, $k\in\Z/(n+3)\Z$. Then, we fill further shifted rows above such
that each diamond
\begin{equation}
  \label{quinconce}
  \begin{matrix}
     & b & \\
    a & & d\\
     & c & 
  \end{matrix}
\end{equation}
verifies $ad=1+bc$, until we reach again a row filled with $1$ only.
\begin{example}\label{frieze}
  Consider the triangulation $T$ of Figure \ref{triangT}.
  \begin{figure}
    \begin{center}
      \scalebox{0.75}{\includegraphics{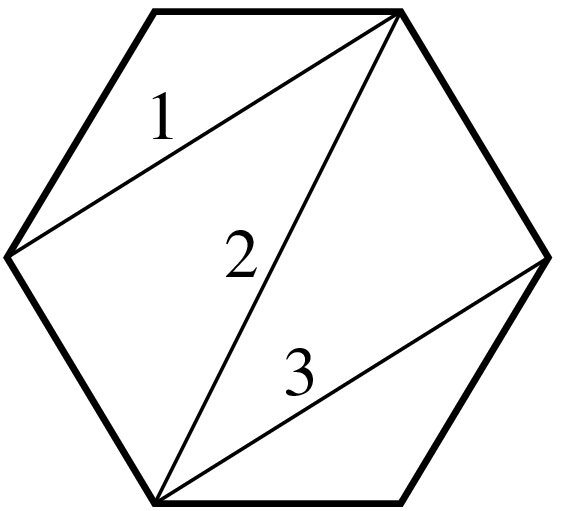}}
      \caption{}
      \label{triangT}
    \end{center}
  \end{figure}
  Then, the corresponding frieze is
\begin{verbatim}
         1     1     1     1     1     1     1  
            1     3     2     1     3     2     1  
         1     2     5     1     2     5     1  
            1     3     2     1     3     2     1  
         1     1     1     1     1     1     1   . 
\end{verbatim}
\end{example}
Now, we make a connection with another construction. In \cite{CCS},
the authors define for each triangulation $T$ of the $(n+3)$-gon a
quiver $Q_T$ in the following way. Let $D_i$, $1\leq i\leq n$, be the
diagonals of $T$. The set of vertices of the quiver $Q_T$ is
$\{1,\ldots,n\}$, and $i\rightarrow j$ if and only if $D_i$ and $D_j$
are edges of a triangle in $T$ and if the angle from $D_i$ to $D_j$ is
counterclockwise. We can choose the triangulation such that $Q_T$ is
any orientation of the Dynkin diagram of type A$_n$. In this case, we
can define the categories $\mkq$ and $\cc$ as before.\par
For any indecomposable object $M$ of $\cc$, let $\mathsf{x}_M$ be the
following number
\begin{equation}
  \mathsf{x}_M=X_M\mid_{u_1=\dots=u_n=1}.
\end{equation}
We have the proposition:
\begin{proposition}
  Let $T$ be a triangulation of the $(n+3)$-gon such that the associated
  quiver $Q_T$ is an orientation of the Dynkin diagram of type A$_n$.
  Let $\gamma_T$ be obtained from the AR-quiver of $\Gamma_{Q_T}$ by
  replacing each $M$ by the number $\mathsf{x}_M$. Then, $\gamma_T$ is
  the Coxeter-Conway frieze associated to $T$.
\end{proposition}
\begin{remark}
  The reader can check the proposition on an example by comparing
  Example \ref{A3} and Example \ref{frieze}.
\end{remark}
  \begin{figure}
    \begin{center}
      \scalebox{0.75}{\includegraphics{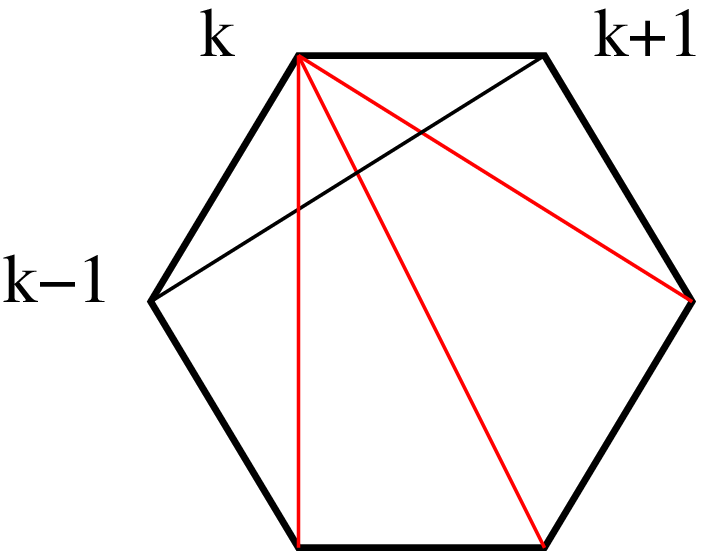}}
      \caption{}
      \label{derniere}
    \end{center}
  \end{figure}
\begin{proof}
  We sketch the proof of the proposition. In \cite[par. 5]{CCS}, the
  authors define a bijection from the set of diagonals of the
  $(n+3)$-gon to the set of indecomposable objects of $\cc$. And
  moreover, the objects of the first row of the AR-quiver $\gcc$ of
  $\cc$ correspond to diagonals of type $[k-1,k+1]$, $k\in\Z/(n+3)\Z$.
  Let $M_k$ be the object corresponding to $[k-1,k+1]$. As $T$ is a
  triangulation, two cases can occur.  Either, $M_k$ is a diagonal of
  $T$, or $M_k$ intersect $T$ non trivially. By \cite[par.  5]{CCS}
  $M_k$ is an object of $\mkq$ if and only if we are in the second
  case. Moreover, in this case the quiver associated to $M_k$ as in
  Example \ref{An} is equioriented. Indeed, the diagonals of $T$
  cutting $[k-1,k+1]$ can be totally ordered in a counterclockwise
  way, as shown in Figure \ref{derniere}. So, by \ref{An}, $M_k$ has
  exactly dim$M_k+1=\mathsf{x}_{M_k}$ submodules. Moreover, the
  dimension of $M_k$ is $[\oplus_i P_i,M_k]=\sum_i[M_k,SP_i]^1$ which,
  by \cite{CCS}, is exactly the number of diagonals of $T$ intersected
  by the corresponding diagonal $[k-1,k+1]$. So,
  $\mathsf{x}_{M_k}=m_k$.  Now, the diamond relations correspond to
  exchange relations as in Proposition \ref{relations}. We obtain the
  result by induction on the rows.
\end{proof}

\bibliographystyle{alpha}
\bibliography{hall}

\end{document}